\documentclass{amsart}
\usepackage{amssymb}
\usepackage{stmaryrd}
\usepackage{fullpage}
\usepackage[colorlinks]{hyperref}

\newtheorem{theorem}{Theorem}[section]
\newtheorem{lemma}[theorem]{Lemma}

\theoremstyle{definition}

\theoremstyle{remark}

\numberwithin{equation}{section}

\newcommand{\bC}{{\mathbf{C}}}

\newcommand{\bF}{{\mathbf{F}}}
\newcommand{\bN}{{\mathbf{N}}}

\newcommand{\Irr}{{\operatorname{Irr}}}

\newcommand{\GL}{{\operatorname{GL}}}
\newcommand{\Gal}{{\operatorname{Gal}}}

\newcommand{\GF}{\mbox{GF}}

\begin{document}

\title{Large orbits of Hall subgroups of solvable linear groups}


\author{Samarth Das}
\address{Heritage High School, Frisco, TX, USA}
\email{samarthdas870@gmail.com}

\author{Yong Yang}
\address{Department of Mathematics, Texas State University, 601 University Drive, San Marcos, TX 78666, USA.}
\email{yang@txstate.edu}

\Large

\date{}



\begin{abstract}
Suppose that $G$ is a finite solvable group and let $H$ be a Hall $\pi$-subgroup, let $b(H)$ be the largest character degree of $H$, we show that $|G:O_{\pi' \pi}(G)|_{\pi} \leq b(H)^2$.
\end{abstract}

\maketitle
\section{Introduction} \label{sec:introduction8}

   Let $G$ be a finite group and $V$ a finite, faithful and completely reducible $G$-module. It is a classical theme to study orbit structures of $G$ acting on $V$. One of the most important and natural questions about orbit structure is to establish the existence of an orbit of a certain size. For a long time, there has been a deep interest and need to examine the size of the largest possible orbits in linear group actions. Using results about length of orbits of linear group action, the following results were proved in ~\cite[Theorem 4.2]{YY7} and ~\cite[Theorem 4.5]{YY7}.

   \begin{theorem}
  Suppose that $G$ is a finite solvable group and let $H$ be a nilpotent Hall $\pi$-subgroup of $G$, then $|G:O_{\pi' \pi}(G)|_{\pi} \leq b(H)^2$.
  \end{theorem}

   \begin{theorem}
  Suppose that $G$ is a finite $\pi$-solvable group and let $H$ be a nilpotent Hall $\pi$-subgroup of $G$, then $|G:O_{\pi' \pi}(G)|_{\pi} \leq b(H)^2$.
  \end{theorem}

  In this paper, we show that the condition of $H$ to be nilpotent could be dropped. We prove the following results.

    \begin{theorem}
  Suppose that $G$ is a finite solvable group and let $H$ be a Hall $\pi$-subgroup of $G$, then $|G:O_{\pi' \pi}(G)|_{\pi} \leq b(H)^2$.
  \end{theorem}

  \begin{theorem}
  Suppose that $G$ is a finite $\pi$-solvable group and let $H$ be a Hall $\pi$-subgroup of $G$, then $|G:O_{\pi' \pi}(G)|_{\pi} \leq b(H)^2$.
  \end{theorem}

If $V$ is a finite vector space of dimension $n$ over $\GF(q)$, where $q$ is a prime power, we denote by $\Gamma(q^n)=\Gamma(V)$ the semilinear group of $V$, i.e.,
\[\Gamma(q^n)=\{x \mapsto ax^{\sigma}\ |\ x \in \GF(q^n), a \in \GF(q^n)^{\times}, \sigma \in \Gal(\GF(q^n)/\GF(q))\}.\]

\section{Main Results} \label{sec:main}

 In order to prove the previously stated results, we will generalize an orbit theorem about Hall $\pi$-subgroups of solvable linear groups. We note that the recent developments in ~\cite{Holt} and ~\cite{Dey} make the calculation relatively easier.

\begin{theorem} \label{thm1}
  Let $G$ be a nontrivial solvable group and let $H$ be a $\pi$-subgroup of $G$. Let $V$ be a faithful $G$-module, over possibly different finite fields of $\pi$-characteristic. Assume that $V_{O_{\pi}(G)}$ is completely reducible. Suppose that $3 \in \pi$, then there exist five orbits of the action of $G$ on $V \oplus V$ with representatives $v_{i_1}+v_{i_2} \in V \oplus V$, $1 \leq i \leq 5$, such that $\bC_H(v_{i_1}) \cap \bC_H(v_{i_2}) \leq O_{\pi}(G)$. Suppose that $3 \not\in \pi$, then there exist three orbits of the action of $G$ on $V \oplus V$ with representatives $v_{i_1}+v_{i_2} \in V \oplus V$, $1 \leq i \leq 5$, such that $\bC_H(v_{i_1}) \cap \bC_H(v_{i_2}) \leq O_{\pi}(G)$.
\end{theorem}
\begin{proof}

If $H \subseteq O_{\pi}(G)$ there is nothing to prove. Thus we may assume that $H \not\subseteq O_{\pi}(G)$. Let $G$ be a counterexample minimizing $\dim(V)$.

Step 1. $V$ is a completely reducible $G$-module.
  Let $R$ be a Hall $\pi'$-subgroup of $O_{\pi \pi'}(G)$. If $h \in HO_{\pi}(G)-O_{\pi}(G)$, let $1 \neq Y(h)$ a Hall $\pi'$-subgroup of $[h,R]$.

  We claim that there exists an irreducible $G$-submodule $V(h)$ of $V$ such that $Y(h)$ acts nontrivially on $V(h)$.

  Since $V_{O_{\pi}(G)}$ is completely reducible and the fields have $\pi$-characteristic, we know that $V_{O_{\pi \pi'}(G)}$ is completely reducible.

  Write $V_{O_{\pi \pi'}(G)}=V_1 \oplus V_2 \cdots \oplus V_t$, where the $V_i$'s are the homogeneous components.

  Since $Y(h)>1$, suppose for instance that $Y(h)$ acts nontrivially on $V_1$. Now consider the $G$-module $\sum_{x \in G}V_1x$ and choose an irreducible $G$-submodule $W$ of it. Let $X$ be an irreducible $O_{\pi \pi'}(G)$-submodule of $W$. Since for every $x \in G$, the $V_1x$'s are homogeneous components, it follows that $X \subseteq V_1x$, for some $x \in G$. Since $Wx=W$, we will have that $V_1 \cap W>0$.

  Suppose now that $Y(h)$ acts trivially on $W \cap V_1$ and let $Y$ be an irreducible $O_{\pi \pi'}(G)$-submodule of $W \cap V_1$. Therefore, since $V_1$ is a direct sum of modules isomorphic to $Y$, it follows that $Y(h)$ acts trivially on $V_1$. A contradiction. This shows that $Y(h)$ acts nontrivially on $V(h)=W$, as claimed.

  Let $U= \sum_{h \in HO_{\pi}(G)-O_{\pi}(G)}V(h)$, $U$ is a completely reducible $G$-module of $V$. If $U<V$, by induction, there exist five orbits or three orbits of the action of $\bar H$ on $U \oplus U$ with representatives $u_{i_1}+u_{i_2} \in U \oplus U$, $1 \leq i \leq 5$ or $3$, such that $\bC_{\bar H}(v_{i_1}) \cap \bC_{\bar H}(v_{i_2}) \leq O_{\pi}(\bar G)$ where $\bar{G}=G/\bC_G(U)$.


  Let $C=\bC_G(U)$ and let $K/C=O_{\pi}(\bar{G})$. Observe that $[K/CO_{\pi}(G),O_{\pi \pi'}(G)/CO_{\pi}(G)]=1$. If $h \in H \cap K-O_{\pi}(G)$, then $[h,R] \subseteq [K,O_{\pi \pi'}(G)] \leq O_{\pi}(G)C$. Since $C$ contains the $\pi'$-subgroups of $O_{\pi}(G)C$, it follows that $Y(h) \subseteq C$, which is a contradiction. This proves that $\bC_{H}(v_{i_1}) \cap \bC_{H}(v_{i_2}) \leq O_{\pi}(\bar G) \leq H \cap K \leq O_{\pi}(G)$ and we may assume $U=V$. Hence $V$ is a completely reducible $G$ module. We note that the case when $3 \not\in \pi$ follows from ~\cite[Theorem 3.2]{YY7} and from now on we may assume that $3 \in \pi$.

Step 2.  $V$ is an irreducible $G$-module with field of characteristic $p$. Assume not, we have $V=V_1 \oplus V_2$ and each $V_i$ is a non-trivial $G$-module. Let $K_i=\bC_G(V_i)$ and $V_i$ is a faithful $G/K_i$-module. By induction, let $v_{i_1}+v_{i_2} \in V_i \oplus V_i$ such that $\bC_{HK_i/K_i}(v_{i_1}) \cap \bC_{HK_i/K_i}(v_{i_2}) \leq O_{\pi}(G/K_i)$ and consider $v_1=v_{1_1}+v_{2_1}$, $v_2=v_{1_2}+v_{2_2}$. Then $\bC_H(v_1) \cap \bC_H(v_2) \leq O_{\pi}(G)$ and the result follows.

Step 3. We now assume that $V$ is not primitive. 

Assume that there exists a proper subgroup $L_1$ of $G$ and an irreducible $L_1$-submodule $V_1$ of $V$ such that $V ={V_1}^G$. We choose $L_1$ to be a maximal subgroup of $G$. Then, $S \cong G/N$ is a primitive permutation group on a right transversal $\Omega$ of $L_1$ in $G$, where $N$ is the normal core of $L_1$ in $G$. Let $V_N=V_1 \oplus \cdots \oplus V_n$, where the $V_i$s are irreducible $L_i$-modules where $L_i=\bN_G(V_i)$ and $m>1$. We know $S$ primitively permutes the $\Omega=\{V_1, \dots, V_n\}$. $K_i=L_i/\bC_G(V_i)$ acts faithfully and irreducibly on $V_i$. We also know that $G$ is isomorphic to a subgroup of $K_1 \wr S$. Let $H$ be a subgroup of $G$ and let $J_i=\bN_H(V_i)/\bC_H(V_i)$.

By induction, $J_i$ has at least $5$ orbits on $V_i \oplus V_i$ for $1 \leq i \leq n$, with representatives $v_{i_{j_1}}+v_{i_{j_2}} \in V_i \oplus V_i$, $1 \leq j \leq 5$, such that $\bC_{J_i}(v_{i_{j_1}}) \cap \bC_{J_i}(v_{i_{j_2}}) \leq O_{\pi}(K_i)$. By mimicking the proof of ~\cite[Proposition 3.2]{YY1}, we can show that $H$ has at least $5$ orbits on $V \oplus V$ with representatives $v_{k_1}+v_{k_2} \in V \oplus V$, $1 \leq k \leq 5$, such that $\bC_H(v_{k_1}) \cap \bC_H(v_{k_2}) \leq O_{\pi}(G)$.

Step 4. We may assume that the action is irreducible and quasi-primitive. We use the notation in \cite[Theorem 2.2]{YY3}. By the main results in ~\cite{Holt} and ~\cite{Dey}, when $e>1$, there are only finite amount of cases left, and we checked those cases using GAP ~\cite{GAP}. By GAP calculation, we note that $H$ has at least $5$ orbits on $V \oplus V$ with representatives $v_{i_1}+v_{i_2} \in V \oplus V$, $1 \leq i \leq 5$, such that $\bC_H(v_{i_1}) \cap \bC_H(v_{i_2}) \leq O_{\pi}(G)$. Note that for all those exceptional cases, we have checked all the possible $\pi$-subgroups to see if it has five orbits on $V \oplus V$ that satisfy the condition. There is only one exception, namely when $G\cong \GL(2,3)$ and $\pi=\{2\}$, and $H$ has exactly three orbits on $V \oplus V$ that satisfy the condition.

Now we may assume $e=1$. By ~\cite[Proposition 2.6]{YY7}, we may assume $p=2$ or $p=3$.

Let $e=1$ and $p=2$. Since $e=1$ we have $G \leq \Gamma(V)=\Gamma(2^d) \cong G_1$ by ~\cite[Corollary 2.3(b)]{manz/wolf}. For any $0 \neq v \in V$, $|\bC_{G_1}(v)| = d$. We can hence assume that $\bC_{G_1}(v)$ is the Galois group of $V = \GF(2^d)$. So the elements of $V$ that do not belong to a regular orbit of $\bC_{G_1}(v)$ are in the union of the subfields $\GF(2^{d/m})$, $m$ varying among the prime divisors of $d$. Since the number of distinct prime divisors of $d$ is at most $\log_2(d)$, it is enough to prove that $f(d) = (2^d-1)-\log_2(d) \cdot (2^{d/2}-1)-4d$ is positive. It is not hard to check that $f(d) > 0$ for all $d \geq 5$. Thus we are left with the cases when $d = 2, 3, 4$ (note $d$ cannot be 1 since the action of $G$ is irreducible). The result can be checked by direct calculations.

Let $e=1$ and $p=3$. Since $e=1$ we have $G \leq \Gamma(V)=\Gamma(3^d) \cong G_1$ by ~\cite[Corollary 2.3(b)]{manz/wolf}. For any $0 \neq v \in V$, $|\bC_{G_1}(v)| = d$. We can hence assume that $\bC_{G_1}(v)$ is the Galois group of $V = \GF(3^d)$. So the elements of $V$ that do not belong to a regular orbit of $\bC_{G_1}(v)$ are in the union of the subfields $\GF(3^{d/m})$, $m$ varying among the prime divisors of $d$. Since the number of distinct prime divisors of $d$ is at most $\log_2(d)$, it is enough to prove that $f(d) = (3^d-1)-\log_2(d) \cdot (3^{d/2}-1)-4d$ is positive. It is not hard to check that $f(d) > 0$ for all $d \geq 3$. Thus we are left with the cases when $d = 1, 2$. The result can be checked by direct calculations.
\end{proof}


\begin{theorem} \label{thm2}
  Suppose that $G$ is a finite $\pi$-solvable group where $\pi$ is a set of primes and let $H$ be a $\pi$-subgroup of $G$. Let $V$ be a faithful $G$-module, over possibly different finite fields of $\pi$-characteristic. Assume that $V_{O_{\pi}(G)}$ is completely reducible, then there exist $v_1,v_2 \in V$ such that $\bC_H(v_1) \cap \bC_H(v_2) \leq O_{\pi}(G)$.
\end{theorem}
\begin{proof}
  If $G$ is solvable, then this is done by Theorem ~\ref{thm1}. It is clear that we may assume $O_{\pi}(G) \subseteq H$ and we denote $\bar H =H/O_{\pi}(G)$. Now let $\bar G =G/O_{\pi}(G)$ and $\bar N=O_{\pi'}(\bar G)$ where $N$ is the preimage of $\bar N$ in $G$. Note that $\bC_{\bar H}(\bar N)=1$ since $\bC_{\bar G}(\bar N) \leq \bar N$. By ~\cite[Theorem 1.2]{MONA}, there exists a nilpotent $\bar H$ invariant subgroup $\bar K$ of $\bar N$ such that $\bC_{\bar H}(\bar K)=1$. Thus we have $O_{\pi}(KH)=O_{\pi}(G)$ and since $V_{KH}$ is faithful, we may assume $G=KH$. Then $G$ is solvable and we are done.
\end{proof}

\begin{lemma} \label{lemeasy}
Suppose that $G$ is a finite group and $V$ is a faithful $G$-module. Assume $G$ has a regular orbit on $V \oplus V$, then there exists $v \in V$ such that $|\bC_G(v)| \leq \sqrt{|G|}$.
\end{lemma}
\begin{proof}
Since there is an element $(v,u) \in V \oplus V$ such that $\bC_G((v, u))= \bC_G(v) \cap \bC_G(u)=1$. Since $|\bC_G(v)| \cdot |\bC_G(u)| = \frac {|\bC_G(v)| \cdot |\bC_G(u)|} {|\bC_G(v) \cap \bC_G(u)|}=|\bC_G(v)  \bC_G(u)| \leq |G|$. It follows that, either  $|\bC_G(v)| \leq \sqrt{|G|}$ or $|\bC_G(u)| \leq \sqrt{|G|}$.
\end{proof}

\begin{theorem} \label{app1}
  Suppose that $G$ is a finite solvable group and let $H$ be a Hall $\pi$-subgroup of $G$, then $|G:O_{\pi' \pi}(G)|_{\pi} \leq b(H)^2$.
\end{theorem}
\begin{proof}
We may assume that $O_{\pi'}(G)=1$. Let $N=O_{\pi}(G)$. Then, fairly standard arguments show that $C=\bC_G(\bF(N)/\Phi(N))\subseteq N$. Write $V=\Irr(\bF(N)/\Phi(N))$ and $\bar{G}=G/C$. Thus $O_{\pi}(\bar{G})=N/C$. Now, $V$ is a faithful $\bar{G}$-module such that $V_{O_{\pi}(\bar{G})}$ is completely reducible.

Let $H$ be a Hall $\pi$-subgroup of $G$ and let $\bar{H}=H/N$. By Theorem ~\ref{thm1} and Lemma ~\ref{lemeasy}, there exists $\lambda \in V$ such that $|\bC_{\bar{H}}(\lambda)| \leq {|\bar{H}|}^{1/2}$. Let $\xi \in \Irr(\bC_{H}(\lambda)|\lambda)$ and $\alpha=\xi^H \in \Irr(H)$. Thus $|G:N|_{\pi} \leq |\bar{H}| \leq \alpha(1)^2 \leq b(H)^2$, as wanted.
\end{proof}

The following result generalizes ~\cite[Theorem 3.2(2)]{MONA}. 

\begin{theorem} \label{app2}
  Suppose that $G$ is a finite $\pi$-solvable group and let $H$ be a Hall $\pi$-subgroup of $G$, then $|G:O_{\pi' \pi}(G)|_{\pi} \leq b(H)^2$.
\end{theorem}
\begin{proof}
The proof is similar to Theorem ~\ref{app1} but using Theorem ~\ref{thm2} instead of Theorem ~\ref{thm1}.
\end{proof}

\section{Acknowledgement} \label{sec:Acknowledgement}
This project was partially supported by a grant from the Simons Foundation (No. 918096, to YY).\\ 

\noindent \textbf{Data availability Statement:} Data sharing not applicable to this article as no datasets were generated or analysed during the current study.\\

\noindent \textbf{Competing interests:} The authors declare none.



\end{document}